\documentclass[11pt]{article}
\usepackage{amsmath,amsthm,amsfonts}
\DeclareMathOperator{\1}{id}
\newcommand\pd[2]{\frac{\partial#1}{\partial#2}}
\renewcommand{\=}{\doteq}
\newcommand{\p}{\partial}
\newcommand{\X}{\mathfrak{X}}
\newcommand{\T}{\mathfrak{T}}
\newtheorem{thm}{Theorem}[section]
\theoremstyle{definition}

\theoremstyle{definition}
\newtheorem{exam}[thm]{Example}
\theoremstyle{definition}
 
\numberwithin{equation}{section}
\numberwithin{equation}{section}
\begin{document}
\title{\bf  Moufang symmetry X.\\
Generalized Lie and Maurer-Cartan equations of\\
continuous Moufang transformations}
\author{Eugen Paal}
\date{}
\maketitle
\thispagestyle{empty}
\begin{abstract}
The differential equations for a continuous birepresentation of a local analytic Moufang loop are established. The commutation relations for the infinitesimal operators of the birepresentation are found. These commutation relations can be seen as a (minimal) generalization of the Maurer-Cartan equations and do not depend on the particular birepresentation.
\par\smallskip
{\bf 2000 MSC:} 20N05, 17D10
\end{abstract}

\section{Introduction}

In this paper we proceed explaing the Moufang symmetry. 
The differential equations for a continuous birepresentation of a local analytic Moufang loop are established. The commutation relations for the infinitesimal operators of the representation are found. These commutation relations can be seen as a (minimal) generalization of the Maurer-Cartan equations and do not depend on the particular birepresentation.

The paper can be seen as a continuation of \cite{Paal1,Paal7}. 

\section{Moufang loops}

A \emph{Moufang loop} \cite{RM} (see also \cite{Bruck,Bel,HP}) is a set $G$ with a binary operation  (multiplication) 
$\cdot: G\times G\to G$,
denoted also by juxtaposition, so that the following three axioms are satisfied:
\begin{enumerate}
\itemsep-3pt
\item[1)] 
in the equation $gh=k$, the knowledge of any two of $g,h,k\in G$ specifies the third one \emph{uniquely},
\item[2)] 
there is a distinguished element $e\in G$ with the property $eg=ge=g$ for all $g\in G$,
\item[3)] 
the \emph{Moufang identity} 
\begin{equation*}
(gh)(kg) = g(hk\cdot g)
\end{equation*}
hold in $G$.
\end{enumerate}
Recall that a set with a binary operation is called a \emph{groupoid}. A groupoid $G$ with axiom 1) is called a \emph{quasigroup}. If axioms 1) and 2) are satisfied, the gruppoid (quasigroup) $G$ is called a \emph{loop}. The element $e$ in axiom 2) is called the \emph{unit} (element) of the (Moufang) loop $G$.

In a (Moufang) loop, the multiplication need not be neither associative nor commutative. Associative (Moufang) loops are well known and are called \emph{groups}. The \emph{associativity} and \emph{commutativity} laws read, respectively,
\begin{equation*}
g(hk)=(gh)k,\quad gh=hg, \qquad \forall g,h,k\in G
\end{equation*}
The associative commutative (Moufang) loops are called the \emph{Abelian} groups.
The most familiar kind of loops are those with the \emph{associative} law, and these are called 
\emph{groups}. A (Moufang) loop $G$ is called \emph{commutative} if the commutativity law holds in $G$, and (only) the commutative associative (Moufang) loops are said to be \emph{Abelian}.

The most remarkable property of the Moufang loops is their \emph{diassociativity}: in a Moufang loop $G$ every two elements generate an associative subloop (group) \cite{RM}. In particular, from this it follows that
\begin{equation} 
\label{alt}
g\cdot gh=g^{2}h,\quad 
hg\cdot g=hg^{2},\quad 
gh\cdot g=g\cdot hg,\qquad \forall g,h\in G
\end{equation}
The first and second identities in (\ref{alt}) are called the left and right \emph{alternativity}, respectively, and the third one is said to be \emph{flexibility}. Note that these identities follow from the Moufang identities as well. 

The unique solution of the equation $xg=e$ ($gx=e$) is called the left (right) \emph{inverse} element of $g\in G$ and is denoted as $g^{-1}_{R}$ ($g^{-1}_{L}$). It follows from the diassociativity of the Moufang loop that 
\begin{equation*}
g^{-1}_{R}=g^{-1}_{L}\doteq  g^{-1},\qquad \forall g\in G\\
\end{equation*}

\section{Analytic Moufang loops and Mal'tsev algebras}

A Moufang loop $G$ is said to be \emph{analytic} \cite{Mal} if $G$ is a finite dimensional real, analytic manifold so that both the Moufang loop operation $G\times G\to G$: ($g,h$) $\mapsto gh$ and the inversion map $G\to G$: $g\mapsto g^{-1}$ are analytic ones. Dimension of $G$ will be denoted as $\dim G\doteq  r$. The local coordinates of $g\in G$ are denoted (in a fixed chart of the unit element $e\in G$) by $g^{1},\dots,g^{r}$, and the local coordinates of the unit $e$ are supposed to be zero: $e^{i}=0$, $i=1,\dots,r$. One has the evident initial conditions
\begin{equation*}
(ge)^{i}=(eg)^{i}=g^{i},\quad i=1,\dots,r
\end{equation*}
As in the case of the Lie groups \cite{Pontr}, we can use the Taylor expansions
%
\begin{align*}
(gh)^{i}
&=h^{i}+u^{i}_{j}(h)g^{j}+\cdots\\
&=g^{i}+v^{i}_{j}(g)h^{j}+\cdots\\
&=g^{i}+h^{i}+a^{i}_{jk}g^{j}h^{k}+\cdots
\end{align*}
%
to introduce the \emph{auxiliary functions} $u^{i}_{j}$ and $v^{i}_{j}$ and the \emph{structure constants}
\begin{equation*} 
c^{i}_{jk}\doteq  a^{i}_{jk}-a^{i}_{kj}=-c^{i}_{kj}   
\end{equation*}
It follows from axiom 1) of the Moufang loop that 
\begin{equation*}
\det(u^{i}_{j})\ne0,\quad
\det(v^{i}_{j})\ne0
\end{equation*}
The \emph{tangent algebra} of $G$ can be defined similarly to the tangent (Lie) algebra of the Lie group 
\cite{Pontr}. Geometrically, this algebra is the tangent space $T_{e}(G)$ of $G$ at $e$. 
The product of $x,y\in T_{e}(G)$ will be denoted by $[x,y]\in T_{e}(G)$. In coordinate form,
\begin{equation*} 
[x,y]^{i}\doteq  c^{i}_{jk}x^{j}y^{k}
            =-[y,x]^{i},\qquad i=1,\dots,r 
\end{equation*}
The tangent algebra will be denoted by  $\Gamma\doteq \{T_{e}(G),[\cdot,\cdot]\}$. The latter algebra need not be a Lie algebra. In other words, there may be a triple $x,y,z\in T_{e}(G)$, such that the Jacobi identity fails:
\begin{equation*} 
J(x,y,z)\doteq  [x,[y,z]]+[y,[z,x]]+[z,[x,y]]\ne0 
\end{equation*}
Instead, for all $x,y,z\in T_{e}(G)$, we have \cite{Mal} a more general identity
\begin{equation*} 
[[x,y],[z,x]]+[[[x,y],z],x]+[[[y,z],x],x]+[[[z,x],x],y]=0
\end{equation*}
called the \emph{Mal'tsev identity}. The tangent algebra is hence said to be the \emph{Mal'tsev algebra}. The Mal'tsev identity concisely reads \cite{Sagle}
\begin{equation*}
[J(x,y,x),x]=J(x,y,[x,z])
\end{equation*}
from which it can be easily seen that every Lie algebra is a Mal'tsev algebra as well. It has been shown in  \cite{Kuzmin} that every  finite-dimensional real Mal'tsev algebra is the tangent algebra of some analytic Moufang loop.

\section{Moufang transformations}

Let $\X$ be a set and let $\T(\X)$ denote the transformation group of $A$. Elements of $\X$ are called \emph{transformations} of $\X$. Multiplication in $\T(\X)$ is defined as the composition of trannsformations, and the unit element of $\T(\X)$ coincides with the identity transformation $\1$ of $A$.

Let $G$ be a Moufang loop with the unit element $e\in G$ and let $(S,T)$ denote a pair of maps $S,T:G\to\T(\X)$. The pair $(S,T)$ is said \cite{Paal7} to be an \emph{action} of $G$ on $\X$ if
\begin{subequations}
\label{bir-def}
\begin{gather}
S_e=T_e=\1\\
S_gT_gS_h=S_{gh}T_g\\
S_gT_gT_h=T_{hg}S_g
\end{gather}
\end{subequations}
hold for all $g,h$ in $G$. The pair $(S,T)$ is called also a \emph{representation} of $G$ (in $\T(\X)$). The transformations $S_g,T_g\in\T(\X)$ ($g\in G$) are called $G$-transformations or the \emph{Moufang transformations} of $\X$.

\begin{exam}
Define the left ($L$) and right ($R$) translations of $G$ by $gh=L_gh=R_hg$.
Then the pair $(L,R)$ of maps $L_g,R_g:G\to \T(G)$ is a representation of $G$.
\end{exam}

Algebraic properties od the Moufang transformations were studied in \cite{Paal7}. In particular, the defining relations
(\ref{bir-def}b,c) can be rewritten as follows:
\begin{equation}
\label{bir-def2}
S_hT_gS_g=T_gS_{hg},\quad
T_hT_gS_g=S_gT_{gh}
\end{equation}
The birepresentation $(S,T)$ is said to be \emph{associative}, if for all $g,h$ in $G$ we have
\begin{equation*}
S_gS_h=S_{gh},\quad
T_gT_h=T_{hg},\quad
S_gT_h=T_hS_g
\end{equation*}
These conditions turn out to be equivalent \cite{Paal7}.

\section{Continuous birepresentations}

Let $G$ be a local analytic Moufang loop and let $\X$ denote a real analytic manifold. We denote dimensions as 
$\dim G\=r$ and $\dim\X\=n$.

An action $(S,T)$ of $G$ on $\X$ is sad to be \emph{differentiable} (smooth, analytic) if the local coordinates of the points $S_gA$ and $T_gA$ are differentiable (smooth, analytic) functions of the points $g\in G$ and $A\in\X$. In this case, the birepresentation is said to be differentiable (smooth, analytic) as well.

In this paper, we consider continuous Moufang transformations only locally, and by 'continuity' we mean differentiability as many times as needed. The action of $g$ from vicinity of the unit element $e$ on $\X$ can be written in local coordinates as
\begin{align*}
(S_gA)^{\mu}
&=S^{\mu}(A^1,\ldots,A^n;g^1,\ldots,g^r)\\
&=S^\mu(A;g)\\
(A_gA)^{\mu}
&=T^{\mu}(A^1,\ldots,A^n;g^1,\ldots,g^r)\\
&=T^\mu(A;g)
\end{align*}
As in case of the Lie transformation groups \cite{Pontr}, we can use the Taylor expansions
\begin{align*}
(S_gA)^{\mu}
&=A^\mu+S^\mu_j(A)g^j+\frac{1}{2}\tilde{S}^\mu_{jk}(A)g^jg^k+O(g^3)\\
(T_gA)^{\mu}
&=A^\mu+T^\mu_j(A)g^j+\frac{1}{2}\tilde{T}^\mu_{jk}(A)g^jg^k+O(g^3)
\end{align*}
to introduce the auxiliary functions $S^\mu_j$ and $T^\mu_j$. The further coefficient in these expansions are assumed to be symmetric with respect to the lower indices:
\begin{equation}
\label{bir-symm}
\tilde{S}^\mu_{jk}=\tilde{S}^\mu_{kj},\quad \tilde{T}^\mu_{jk}=\tilde{T}^\mu_{kj},\quad \text{etc}
\end{equation}

\section{Associators}

An action of $G$ need not be associative even in case $G$ is a group. Nonassociativity of $(S,T)$ can be measured by the formal functions
\begin{align*}
l^\mu(A;g,h)
&\=(S_{gh}A)^{\mu}-(S_gS_hA)^\mu\\
r^\mu(A;g,h)
&\=(T_{gh}A)^{\mu}-(S_hS_gA)^\mu\\
m^\mu(A;g,h)
&\=(T_hS_gA)^{\mu}-(S_gT_hA)^\mu
\end{align*}
which are called associators of $(S,T)$. We have the evident initial conditions
\begin{align*}
l^\mu(A;e,g)&=r^\mu(A;e,g)=m^\mu(A;e,g)\\
l^\mu(A;g,e)&=r^\mu(A;g,e)=m^\mu(A;g,e)
\end{align*}
The associators of $(S,T)$ are considered as the generating expressions in the following sense. First define the 
\emph{first-order} associators $l^{\mu}_{j}$, $\hat{l}^{\mu}_{j}$, $r^{\mu}_{j}$, $\hat{r}^{\mu}_{j}$, $m^{\mu}_{j}$, $\hat{m}^{\mu}_{j}$  by
\begin{align*}
l^{\mu}(A;g,h)
&\doteq  l^{\mu}_{j}(A;h)g^{j}+O(g^{2})\\
&\doteq  \hat{l}^{\mu}_{j}(A;g)h^{j}+O(h^{2})\\
r^{\mu}(A;g,h)
&\doteq  r^{\mu}_{j}(A;h)g^{j}+O(g^{2})\\
&\doteq  \hat{r}^{\mu}_{j}(A;g)g^{j}+O(h^{2})\\
m^{\mu}(A;g,h)
&\doteq  m^{\mu}_{j}(A;h)g^{j}+O(g^{2})\\
&\doteq  \hat{m}^{\mu}_{j}(A,g)h^{j}+O(h^{2})
\end{align*} 
As an example, calculate $l^{\mu}_{j}$. We have
\begin{align*}
(S_gS_hA)^\mu
&=S^\mu(S_hA;g)\\
&=(S_hA)^\mu+S^{\mu}_{j}(S_hA)j^j+O(g^2)\\
(S_{gh}A)^\mu
&=S^\mu(A;gh)\\
&=(S_{h}A)^{\mu}+\pd{(S_gA)^\mu}{h^k}u^k_j(h)g^j+O(g^2)
\end{align*}
so that
\begin{equation*}
l^{\mu}_{j}(A;h)=u^k_j(h)g^j \pd{(S_gA)^\mu}{h^k} - S^\mu_j(S_gA)
\end{equation*}
Remaining first-order associators can be found similarly and result read
\begin{subequations}
\label{a1}
\begin{align}
l^{\mu}_{j}(A;g) 
& =u^{s}_{j}(g)\pd{(S_gA)^{\mu}}{g^{s}}-S^{\mu}_{j}(S_gA)\\
\hat{l}^{\mu}_{j}(A;g) 
& =v^{s}_{j}(g)\pd{(S_gA)^{\mu}}{g^{s}}-S^{\nu}_{j}(A)\pd{(S_gA)^{\mu}}{A^{\nu}}\\
r^{\mu}_{j}(A;g) 
&=u^{s}_{j}(g)\pd{(T_gA)^{\mu}}{g^{s}}-T^{\nu}_{j}(A)\pd{(T_gA)^{\mu}}{A^{\nu}}\\
\hat{r}^{\mu}_{j}(A;g)
&=v^{s}_{j}(g)\pd{(T_gA)^{\mu}}{g^{s}}-T^{\mu}_{j}(T_gA)\\
m^{\mu}_{j}(A;g)  
&=-S^{\mu}_{j}(T_gA)+S^{\nu}_{j}(S_gA)\pd{(T_gA)^{\mu}}{A^{\nu}}\\
\hat{m}^{\mu}_{j}(A;g) 
&=-T^{\mu}_{j}(S_gA)+T^{\nu}_{j}(S_gA)\pd{(S_gA)^{\mu}}{A^{\nu}}
\end{align}
\end{subequations}
Next one can check the initial conditions
\begin{align*}
l^{\mu}_{j}(A;e)&=r^{\mu}_{j}(A;e)=m^{\mu}_{j}(A;e)=0\\
l^{\mu}_{j}(A;e)&=r^{\mu}_{j}(A;e)=m^{\mu}_{j}(A;e)=0
\end{align*}
and define the \emph{second-order} associators
$l^{\mu}_{jk}$, $\hat l^{\mu}_{jk}$, $m^{\mu}_{jk}$, $\hat m^{\mu}_{jk}$, $r^{\mu}_{jk}$, $\hat r^{\mu}_{jk}$ by
\begin{subequations}
\label{bir-2-assoc}
\begin{align} 
l^{\mu}_{j}(A;g)
&\doteq  l^{\mu}_{jk}(A)g^{k}+O(g^{2})\\
&\doteq  \hat l^{\mu}_{jk}(A)g^{k}+O(g^{2})\\
r^{\mu}_{j}(A;g)
&\doteq  r^{\mu}_{jk}(A)g^{k}+O(g^{2})\\
&\doteq  \hat r^{\mu}_{jk}(A)g^{k}+O(g^{2})\\
m^{\mu}_{j}(A;g)
&\doteq  \hat m^{\mu}_{jk}(A)g^{k}+O(g^{2})\\
 &\doteq  m^{\mu}_{jk}(A)g^{k}+O(g^{2})
\end{align}
\end{subequations}
Calculating, we get
\begin{subequations}
\label{a2}
\begin{align} 
l^{\mu}_{jk}(A)
&=\hat m^{\mu}_{kj}(A)   \\
&=\tilde{S}^{\mu}_{kj}(A)-a^{\nu}_{jk}S^{\mu}_{s}(A) + S^{\nu}_{k}(A)\pd{S^{\mu}_{j}(A)}{A^{\nu}} \\
m^{\mu}_{jk}(A)
&=\hat r^{\mu}_{kj}(A)  \\
&=\tilde{T}^{\mu}_{jk}(A)+a^{\nu}_{jk}T^{\mu}_{s}(A) - T^{\nu}_{j}(A)\pd{T^{\mu}_{k}(A)}{A^{\nu}}  \\
r^{\mu}_{jk}(A)
&=-\hat{m}^{\mu}_{kj}(A)  \\
&= S^{\nu}_{j}(A)\pd{T^{\mu}_{k}(A)}{A^{\nu}} - T^{\nu}_{k}(A)\pd{S^{\mu}_{j}(A)}{A^{\nu}}
\end{align}
\end{subequations}

\section{Minimality conditions and generalized Lie equations}

Differentiate the definig relation (\ref{bir-def}b,c) of a birepresntation $(S,T)$ and relations (\ref{bir-def2})
in local coordinates with respect to $g^{j}$ at $g=e$. Then, redenoting $h\to g$  we obtain for the first-order associators constraints
\begin{subequations}
\label{m-bir-d1}
\begin{align} 
\hat{l}^{\mu}_{j}(A;g)=\hat{m}^{\mu}_{j}(A;g)=-l^{\mu}_{j}(A;g)\\
\hat{r}^{\mu}_{j}(A;g)=\hat{m}^{\mu}_{j}(A;g)=-\hat{r}^{\mu}_{j}(A;g)
\end{align}
\end{subequations}
As an example, check relation $\hat{m}^{\mu}_{j}=-l^{\mu}_{j}$.
Write the defing relation (\ref{bir-def}b) in local coordinates:
\begin{equation*}
(S_gT_gS_hA)^\mu=(S_{gh}T_gA)^\mu,\quad \mu=1,\ldots,n
\end{equation*}
Calculate:
\begin{align*}
(S_gT_gS_hA)^\mu
&=S^\mu(T_gS_hA;g)\\
&=(T_gS_hA)^\mu+S^\mu_j(T_gS_hA)g^j+O(g^2)\\
&=(S_hA)^\mu+T^\mu_j(S_hA)g^j+S^\mu_j(S_hA)g^j+O(g^2)\\
(S_{gh}T_gA)^\mu
&=S^\mu(T_gA;gh)\\
&=(S_hA)^\mu+\pd{(S_hA)^\mu}{A^\nu}T^\nu_j(A)g^j+\pd{(S_hA)^\mu}{h^s}u^s_j(h)g^j+O(g^2)
\end{align*}
Comparing the above expansions we get the desired relation $\hat{m}^{\mu}_{j}=-l^{\mu}_{j}$.
Remaining relations from (\ref{m-bir-d1}a,b) can be checked similarly.

If the birepresentation $(S,T)$ is required to be associative, we get the familiar Lie equations \cite{Pontr} of a Lie transformation group:
\begin{align*} 
\hat{l}^{\mu}_{j}(A;g)=\hat{m}^{\mu}_{j}(A;g)=-l^{\mu}_{j}(A;g)=0\\
\hat{r}^{\mu}_{j}(A;g)=\hat{m}^{\mu}_{j}(A;g)=-\hat{r}^{\mu}_{j}(A;g)=0
\end{align*}
In a sense, one may say that birepresentations of the Moufang loop have the property of the 'minimal' deviation from associativity. Thus the differential identities (\ref{m-bir-d1}a,b) are called the 
\emph{first-order minimality conditions} of $(S,T)$.

The first-order minimality conditions  (\ref{m-bir-d1}a,b) read as the differential equations for $G$-trans\-formations.
Define the auxiliary functions $w^s_j$ and  $P^\mu_j(g)$ by 
\begin{gather}
\label{uvw}
S^s_j(g)+T^s_j(g)+w^s_j(g)=0\\
\label{STP}
S^\mu_j(A)+T^\mu_j(A)+P^\mu_j(A)=0
\end{gather}
For $S_gA$ the \emph{generalized Lie equations} (GLE) read
\begin{subequations}
\label{gle_S}
\begin{align}
u^{s}_{j}(g)\pd{(S_gA)^{\mu}}{A^{\nu}}+T^{\nu}_{j}(A)\pd{(S_gA)^{i}}{A^{\nu}}+P^{\nu}_{j}(S_gA)&=0\\
v^{s}_{j}(g)\pd{(S_gA)^{\mu}}{A^{\nu}}+P^{\nu}_{j}(h)\pd{(S_gA)^{i}}{A^{\nu}}+T^{\nu}_{j}(S_gA)&=0\\
w^{\nu}_{j}(g)\pd{(S_gA)^{\mu}}{A^{\nu}}+S^{\nu}_{j}(h)\pd{(S_gA)^{i}}{A^{\nu}}+S^{\nu}_{j}(S_gA)&=0
\end{align}
\end{subequations}
For $T_gA$ the GLE read
\begin{subequations}
\label{gle_T}
\begin{align}
v^{s}_{j}(g)\pd{(T_gA)^{\mu}}{A^{\nu}}+S^{\nu}_{j}(A)\pd{(T_gA)^{i}}{A^{\nu}}+P^{\nu}_{j}(T_gA)&=0\\
u^{s}_{j}(g)\pd{(T_gA)^{\mu}}{A^{\nu}}+P^{\nu}_{j}(h)\pd{(T_gA)^{i}}{A^{\nu}}+S^{\nu}_{j}(T_gA)&=0\\
w^{\nu}_{j}(g)\pd{(T_gA)^{\mu}}{A^{\nu}}+T^{\nu}_{j}(h)\pd{(T_gA)^{i}}{A^{\nu}}+T^{\nu}_{j}(T_gA)&=0
\end{align}
\end{subequations}
Due to (\ref{uvw}) and  (\ref{STP}) these differential equations are linearly dependent: by adding  (\ref{gle_S}a--c)
or (\ref{gle_T}a--c) we get $0=0$.

\section{Generalized Maurer-Cartan equations}

Differentiate constraints (\ref{m-bir-d1}a--c)  with respect to $g^{k}$  at $g=e$. Then we obtain
\begin{align*} 
\hat{l}^{\mu}_{jk}(A)=\hat{m}^{\mu}_{jk}(A)=-l^{\mu}_{jk}(A)=0\\
r^{\mu}_{jk}(A)=m^{\mu}_{jk}(A)=-\hat{r}^{\mu}_{jk}(A)=0
\end{align*}
Using here (\ref{bir-2-assoc}a,c,e)  we obtain the \emph{second-order minimality conditions} of $(S,T)$:
\begin{equation*}
\hat{l}^{\mu}_{jk}(g)=r^{\mu}_{jk}(g)=m^{\mu}_{jk}(g)=-m^{\mu}_{kj}(g)
\end{equation*}
Again, for associative $G$-transformations we have
\begin{equation*}
\hat{l}^{\mu}_{jk}(g)=r^{\mu}_{jk}(g)=m^{\mu}_{jk}(g)=-m^{\mu}_{kj}(g)=0 
\end{equation*}
which justifies the term 'minimality conditions'.

It follows from skew-symmetry $\hat{l}^{\mu}_{jk}=-\hat{l}^{\mu}_{kj}$ and $r^{\mu}_{jk}=-r^{\mu}_{kj}$, respectively, that
\begin{subequations}
\begin{align}
\label{st-bir}
2\tilde{S}^{\mu}_{jk}&=S^{\nu}_{k}\pd{S^{\mu}_{j}}{A^{\nu}}+S^{\nu}_{j}\pd{S^{\mu}_{k}}{A^{\nu}}
            -\left(a^{s}_{jk}+a^{s}_{kj}\right)S^{\mu}_{s}\\
2\tilde{T}^{\mu}_{jk}&=T^{\nu}_{k}\pd{T^{\mu}_{j}}{A^{\nu}}+T^{\nu}_{j}\pd{T^{\mu}_{k}}{A^{\nu}}
            -\left(a^{s}_{jk}+a^{s}_{kj}\right)T^{\mu}_{s}
\end{align}
\end{subequations}
Note thet here we used  the symmetry  property (\ref{bir-symm}) as well.
Express $\tilde{S}^{\mu}_{jk}$ and $\tilde{T}^{\mu}_{jk}$ from these relations and substitute into (\ref{bir-2-assoc}b)  and (\ref{bir-2-assoc}d), respectively. The result reads
\begin{align*}
S^{\nu}_{k}\pd{S^{\mu}_{j}}{A^{\nu}}-S^{\nu}_{j}\pd{S^{\mu}_{k}}{A^{\nu}}
             &=c^{s}_{jk}T^{\mu}_{s}+2\hat{l}^{\mu}_{jk}\\
T^{\nu}_{k}\pd{T^{\mu}_{j}}{A^{\nu}}-T^{\nu}_{j}\pd{T^{\mu}_{k}}{A^{\nu}}
             &=c^{s}_{kj}T^{\mu}_{s}+2r^{\mu}_{jk}
\end{align*}
Now using the equalities $\hat{l}^{\mu}_{jk}=m^{\mu}_{jk}$, $r^{\mu}_{jk}=-m^{\mu}_{kj}$
and formula (\ref{bir-2-assoc}f) for $m^{\mu}_{jk}$, we obtain the differential equations for 
the auxiliary functions $S^{\mu}_{j}$ and  $T^{\mu}_{j}$:
\begin{subequations}
\label{gen-m-c-bir}
\begin{align}
S^{\nu}_{k}\pd{S^{\mu}_{j}}{A^{\nu}}-S^{\nu}_{j}\pd{S^{\mu}_{k}}{A^{\nu}}
&=c^{s}_{jk}S^{\mu}_{s} +2\left(T^{\nu}_{j}\pd{S^{\mu}_{k}}{A^{\nu}}
                           -S^{\nu}_{k}\pd{T^{\mu}_{j}}{A^{\nu}}\right)\\
T^{\nu}_{k}\pd{T^{\mu}_{j}}{A^{\nu}}-T^{\nu}_{j}\pd{T^{\mu}_{k}}{A^{\nu}}
&=c^{s}_{kj}T^{\mu}_{s}+2\left(S^{\nu}_{j}\pd{T^{\mu}_{k}}{A^{\nu}}
                           -T^{\nu}_{k}\pd{S^{\mu}_{j}}{A^{\nu}}\right)
\end{align}
\end{subequations}
called the generalized Maurer-Cartan equations for $G$-transformations.
In a sense, the generalized Maurer-Cartan equations generalize the Maurer-Cartan equations \cite{Pontr} in the minimal way.

The generalized Maurer-Cartan differential equations can be rewritten more concisely. 
For $x\in T_e(G)$ introduce the \emph{infinitesimal $G$-transformations}:
\begin{equation*}
S_x\doteq x^j S^{\mu}_j(g)\frac{\p}{\p A^{\mu}},\quad 
T_y\doteq x^j T^{\mu}_j(g)\frac{\p}{\p A^{\mu}}\quad \in T_g(G)
\end{equation*}
Then  the generalized Maurer-Cartan  equations (\ref{gen-m-c-bir}a,b) can be rewritten, respectively, as  the commutation relations
\begin{subequations}
\label{s-t-xy}
\begin{align}
\left[L_{x},L_{y}\right]&=L_{[x,y]}-2\left[L_{x},R_{y}\right] \\
\left[R_{x},R_{y}\right]&=R_{[y,x]}-2\left[R_{x},L_{y}\right]\\
\left[L_{x},R_{y}\right]&=\left[R_{x},L_{y}\right],\quad x,y\in T_(G)
\end{align}
\end{subequations}
Note that  commutation relation (\ref{s-t-xy}c) can easily be obtained from the identities
\begin{equation*}
\left[S_{x},S_{y}\right]=-\left[S_{y},S_{x}\right],
\quad
\left[T_{x},T_{y}\right]=-\left[T_{y},T_{x}\right]
\end{equation*}
Thus, finally the (generalized) Maurer-Cartan equations for infinitesimal $G$ transformations read
\begin{equation*}
2\left[S_{x},T_{y}\right]
=S_{[x,y]}-\left[S_{x},S_{y}\right]  
=T_{[y,x]}-\left[T_{x},T_{y}\right]
=2\left[T_{x},S_{y}\right]=0 
\end{equation*}

\section*{Acknowledgement}

Research was in part supported by the Estonian Science Foundation, Grant 6912.

\bigskip\noindent
Department of Mathematics\\
Tallinn University of Technology\\
Ehitajate tee 5, 19086 Tallinn, Estonia\\ 
E-mail: eugen.paal@ttu.ee

\end{document}